\documentstyle{amsppt}
\baselineskip18pt
\magnification=\magstep1
\pagewidth{30pc}
\pageheight{45pc}
\hyphenation{co-deter-min-ant co-deter-min-ants pa-ra-met-rised
pre-print pro-pa-gat-ing pro-pa-gate
fel-low-ship Cox-et-er dis-trib-ut-ive}
\def\leaderfill{\leaders\hbox to 1em{\hss.\hss}\hfill}

\def\H{{\Cal H}}

\def\idest{i.e.,\ }

\def\a{{\alpha}}
\def\be{{\beta}}

\def\l{{\lambda}}

\def\s{{\sigma}}

\def\W{{\widehat W}}

\def\b0{\text{\bf 0}}

\def\ds{\displaystyle}

\def\complex{{\Bbb C}}
\def\zed{{\Bbb Z}}

\def\enn{{\Bbb N}}

\def\boxit#1{\vbox{\hrule\hbox{\vrule \kern3pt
\vbox{\kern3pt\hbox{#1}\kern3pt}\kern3pt\vrule}\hrule}}
\def\rabbit{\vbox{\hbox{\kern0pt
\vbox{\kern0pt{\hbox{---}}\kern3.5pt}}}}

\def\tableau#1{
        \hbox {
                \hskip -10pt plus0pt minus0pt
                \raise\baselineskip\hbox{
                \offinterlineskip
                \hbox{#1}}
                \hskip0.25em
        }
}

\def\tabCol#1{
\hbox{\vtop{\hrule
\halign{\strut\vrule\hskip0.5em##\hskip0.5em\hfill\vrule\cr\lower0pt
\hbox\bgroup$#1$\egroup \cr}
\hrule
} } \hskip -10.5pt plus0pt minus0pt}

\def\CR{
        $\egroup\cr
        \noalign{\hrule}
        \lower0pt\hbox\bgroup$
}



\def\blank#1#2{
\hbox to #1{\hfill \vbox to #2{\vfill}}
}

\def\domeq{\trianglerighteq}


\def\strut{\vrule height10pt depth5pt width0pt}

\topmatter
\title On $321$-avoiding permutations in affine Weyl groups
\endtitle

\author R.M. Green \endauthor
\affil 
Department of Mathematics and Statistics\\ Lancaster University\\
Lancaster LA1 4YF\\ England\\
{\it  E-mail:} r.m.green\@lancaster.ac.uk\\
\endaffil

\abstract
We introduce the notion of $321$-avoiding permutations in the affine
Weyl group $W$ of type $A_{n-1}$ by considering the group as a George
group (in the sense of Eriksson and Eriksson).  This enables us to
generalize a result of Billey, Jockusch and Stanley to show that the
$321$-avoiding permutations in $W$ coincide with the
set of fully commutative elements; in other words, any two reduced
expressions for a $321$-avoiding element of $W$ (considered as a
Coxeter group) may be obtained from each other by repeated applications
of short braid relations.  

Using Shi's characterization of the
Kazhdan--Lusztig cells in the group $W$, we use our main result to
show that the fully
commutative elements of $W$ form a union of Kazhdan--Lusztig cells.
This phenomenon has been studied by the author and J. Losonczy for
finite Coxeter groups, and is interesting partly because it allows
certain structure constants for the Kazhdan--Lusztig basis of the
associated Hecke algebra to be computed combinatorially.  

We also show how some of our results can be generalized to a 
larger group of permutations, the extended affine Weyl group
associated to $GL_n(\complex)$.
\endabstract

\thanks
This paper was written at the Newton Institute, Cambridge,
and supported by the special semester on Symmetric Functions and
Macdonald Polynomials, 2001.
\endthanks

\endtopmatter

\centerline{\bf To appear in the Journal of Algebraic Combinatorics}

\head Introduction \endhead

Let $W$ be a Coxeter group with generating set 
$S = \{s_i\}_{i \in I}$.  The fully commutative elements of $W$ were defined by
Stembridge, and may be characterised \cite{{\bf 18}, Proposition 1.1} as 
those elements all of whose reduced
expressions avoid strings of length $m$ of the form $s_i s_j s_i
\cdots$, where $2 < m = m_{i, j}$ is the order of $s_i s_j$.
Billey, Jockusch and Stanley \cite{{\bf 1}, Theorem 2.1} show that in the
case where $W$ is the symmetric group (under the natural
identifications), $w \in W$ is fully commutative
if and only if it is $321$-avoiding.  The latter condition means
that there are no $i < j < k$ such that $w(i) > w(j) > w(k)$.   It
is clear from properties of root systems of type $A$ that the
$321$-avoiding elements $w \in W$ are precisely those for which
there are no positive roots $\a, \be$ such that both (i) $\a + \be$ is a root
and (ii) $w(\a)$ and $w(\be)$ are both negative roots.

In this paper we consider the affine analogue of this problem.  To
define the analogue of a $321$-avoiding permutation, we recall
Lusztig's realization of the affine Weyl group of type
$\widetilde A_{n-1}$ as a group
of permutations of the integers.  The interplay between the
permutation representation and the root system of type $\widetilde A$ can
be seen in the work of Papi \cite{{\bf 15}}, and the relationship
between full commutativity and root systems in simply laced types has
been established by Fan and Stembridge \cite{{\bf 6}}.  It is therefore not
hard to show that the $321$-avoiding elements and the fully
commutative elements coincide in type $\widetilde A$, but we choose to
give a different derivation of this result (avoiding root systems) 
because the equivalence of 
conditions (iii) and (iv) of Theorem 2.7 is not made explicit in 
\cite{{\bf 15}}.

Our real interest in this problem stems from a desire to understand
the Kazhdan--Lusztig basis combinatorially.
The two-sided Kazhdan--Lusztig cells for
Coxeter groups of type $\widetilde A_{n-1}$ were classified by Shi
\cite{{\bf 16}} according to a combinatorial criterion on elements $w \in W$
generalizing the notion of $321$-avoidance.  Although it may be possible
to classify the Kazhdan--Lusztig cells of a given Coxeter group using
a purely algebraic or combinatorial approach, it is not reasonable to
expect that the Kazhdan--Lusztig basis $\{C'_w : w \in W\}$
itself and its structure constants may
be understood using similar means.  However, as we discuss in \S5, 
if the fully commutative
elements form a union of Kazhdan--Lusztig cells closed under
$\geq_{LR}$, results of the author and J. Losonczy in \cite{{\bf 10}}
reduce the problem of computing the coefficient of $C'_z$ in $C'_x
C'_y$ to a tractable combinatorial problem in the case that $z$ is
fully commutative.

We also show how these results may be generalized to a certain
extension of the affine Weyl group of type $\widetilde A_{n-1}$.

\head 1. Affine Weyl groups of type $A$ as permutations \endhead

In \S1, we show how the affine Weyl group $W$ of type $A$ and its extension
$\W$ may be viewed as groups of permutations of $\zed$.  In the
case of the Coxeter group $W$, Eriksson and Eriksson
\cite{{\bf 3}, \S5} call such a group of permutations a ``George group''
(after G. Lusztig).  The group $\W$ is called (following \cite{{\bf 19},
\S2.1}) the extended affine Weyl group associated to $GL_n(\complex)$,
and it can be similarly considered as a group of permutations.  We
follow the treatment of these groups in \cite{{\bf 7},
\S1} but work with left actions rather than right actions.

The affine Weyl group of type $\widetilde A_{n-1}$
arises from the Dynkin diagram in Figure 1.

\topcaption{Figure 1} Dynkin diagram of type $\widetilde A_{n-1}$ 
\endcaption
\centerline{
\hbox to 2.958in{
\vbox to 0.944in{\vfill
        \includegraphics{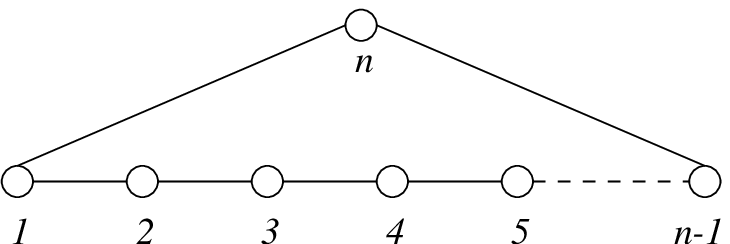}
}
\hfill}
}

The number of vertices in Figure 1 is $n$.  The top
vertex is regarded as an extra relative to the
remainder of the graph, which is a Coxeter graph of type $A_{n-1}$.
We will assume $n \geq 3$ throughout this paper.

We associate an affine Weyl group, $W = W(\widetilde A_{n-1})$,
to this Dynkin diagram in the usual way (as in \cite{{\bf 12}, \S2.1}).  
This associates to node $i$
of the graph a generating involution $s_i$ of $W$, where $s_i s_j =
s_j s_i$ if $i$ and $j$ are not adjacent in the graph, and $$
s_i s_j s_i = s_j s_i s_j
$$ if $i$ and $j$ are adjacent in the graph.  For $t \in {\Bbb Z}$, 
it is convenient to
denote by $\bar{t}$ the congruence class of $t$ modulo $n$, taking
representatives in the set $\{1, 2, \ldots, n\}$.
For the purposes of this paper, it is helpful to think of
the group $W$ as follows.

\proclaim{Proposition 1.1 (Lusztig)}
There exists a group isomorphism from $W$
to the set of permutations of ${\Bbb Z}$ which satisfy
the following conditions: $$
\eqalignno{
w(i+n) &= w(i) + n \text{ for } i \in {\Bbb Z}, & (1)\cr
\sum_{t = 1}^n w(t) &= \sum_{t = 1}^n t & (2)\cr
}$$
such that $s_i$ is mapped to the permutation
$$t \mapsto \cases
t & \text{ if } \bar{t} \ne \bar{i}, \overline{i + 1},\cr
t-1 & \text{ if } \bar{t} = \overline{i + 1},\cr
t+1 & \text{ if } \bar{t} = \bar{i},\cr
\endcases$$ for $t \in {\Bbb Z}$.
\endproclaim

\demo{Proof}
This result is stated in \cite{{\bf 14}, \S3.6}; a proof may be found in
\cite{{\bf 3}, Theorem 20}.
\qed\enddemo

We also define an extension, $\W$, of $W$ as a 
group of permutations of ${\Bbb Z}$, as in \cite{{\bf 19}}.  This group
will be discussed in more depth in \S4.

\definition{Definition 1.2}
Let $\rho$ be the permutation of ${\Bbb Z}$ taking $t$ to $t + 1$ for
all $t$.  Then the group $\W$ is defined to be the group of
permutations of ${\Bbb Z}$ generated by the group $W$ and $\rho$.
\enddefinition

\proclaim{Proposition 1.3}
There exists a group isomorphism from $\W$
to the set of permutations of ${\Bbb Z}$ that satisfies the
condition $w(i+n) = w(i)$ for all $i \in {\Bbb Z}$.
Any element of $\W$ is uniquely expressible in the form $\rho^z w$ for
$z \in {\Bbb Z}$ and $w \in W$.  Conversely, any element of this form
is an element of $\W$.
\endproclaim

\demo{Proof}
Note that we automatically have $$
\sum_{t = 1}^n w(t) \equiv \sum_{t = 1}^n t \mod n
$$ because $w$ is a permutation.  The result now follows from 
\cite{{\bf 7}, Proposition 1.1.3, Corollary 1.1.4}.
\qed\enddemo

The group $\W$ can be expressed in a more familiar way using a
semidirect product construction.

\proclaim{Proposition 1.4}
Let ${\Cal S}_n$ be the subgroup of $\W$ generated by $$
\{s_1, s_2, \ldots, s_{n-1}\}
.$$  Let $Z$ be the
subgroup of $\W$ consisting of all permutations $z$ satisfying $$
z(t) \equiv t \mod n
$$ for all $t$.  Then ${\Bbb Z}^n \cong Z \trianglelefteq \W$ and $W$
is isomorphic to the semidirect product of ${\Cal S_n}$ and $Z$.
\endproclaim

\demo{Proof}
See \cite{{\bf 7}, Proposition 1.1.5}.
\qed\enddemo

In the language of \cite{{\bf 3}, 5.1.2}, the George group $W$ 
is defined as the group of ``locally finite'' permutations which 
commute with all elements of the ``rigid group'' of translations 
generated by $\rho^n$.  Proposition 1.4 is saying that the group $\W$
may be obtained by dropping the hypothesis of local finiteness.

It is convenient to extend the usual notion of the length of an
element of a Coxeter group to the group $\W$ in the following way.

\definition{Definition 1.5}
Let $w \in W$ and let $w = s_{i_1} s_{i_2} \cdots s_{i_r}$ be a word in
the set $S$ of generating involutions such that $r$ is minimal.  Then we
define the length of $w$, $\ell(w)$, to be equal to $r$.

The length, $\ell(w')$, of a typical element $w' = \rho^z w$ of $\W$
(where $z \in {\Bbb Z}$ and $w \in W$) is defined to be $\ell(w)$.
\enddefinition

Definition 1.5 is very natural from the point of view of George
groups: in the language of \cite{{\bf 3}, \S6.1} it states that the length
of an element $w \in \W$ is equal to the ``class
inversion number'' of $w$.

\head 2. Pattern avoidance in Coxeter groups \endhead

In \S2, we restrict our attention to the Coxeter group $W$ of type 
$\widetilde A_{n-1}$.  We first recall the notions of full
commutativity and $321$-avoidance.

\definition{Definition 2.1}
Let $w \in W$.  

We say $w$ is {\it fully commutative} if no reduced
expression for $w$ contains a subword of the form $s_i s_j s_i$
for Coxeter generators $s_i$ and $s_j$.  The set of fully commutative
elements of $W$ is denoted by $W_c$.

We say $w$ is $321$-avoiding if, when $W$ is considered as a George
group as in Proposition 1.1, there are no integers $a < b < c$ such
that $w(a) > w(b) > w(c)$.
\enddefinition

\remark{Remark 2.2}
Notice that if $w \in W$ is such that $a < b < c$ and $w(a) > w(b) >
w(c)$, it must be the case that $w(a)$, $w(b)$ and $w(c)$ are pairwise
non-congruent modulo $n$.
\endremark

The following result will be used to address the issue of checking whether
an infinite permutation is $321$-avoiding.

\proclaim{Proposition 2.3}
Let $w \in W$ and $a, b, c \in \zed$ with $a < b < c$.
If $w(a) > w(b) > w(c)$, then there
exist $a', c' \in \zed$ with $0 < b - a' < n$, $0 < c' - b < n$,
$a \equiv a' \mod n$, $c \equiv c' \mod n$, $a' < b < c'$ and
$w(a') > w(b) > w(c')$.
\endproclaim

\demo{Proof}
First, note that $a$, $b$ and $c$ are pairwise noncongruent modulo $n$ by
Remark 2.2.

Write $c = c' + kn$ with $k \in \zed^{\geq 0}$ and $0 < c' - b < n$.
By condition (1) of Proposition 1.1, we may assume without loss of generality
that $k = 0$ by replacing the triple $(a, b, c)$ by $(a, b, c')$.  This
satisfies the claimed properties because $w(b) > w(c) = w(c-kn) + kn \geq w(c')$.

The argument to find $a'$ is similar, starting with the condition
$a = a' - kn$ with $k \in \zed^{\geq 0}$ and $0 < b - a' < n$.
\qed\enddemo

\remark{Remark 2.4}
Suppose $w \in W$ is given as a permutation.  
Proposition 2.3 shows that we can determine whether there 
exist $a < b < c$ with $w(a) > w(b) > w(c)$ by testing
$\ds{{n \choose 3}}$ triples.  This is because if such a triple exists, we
may assume that $1 \leq b \leq n$ by condition (1) of Proposition 1.1, and 
we may apply
Proposition 2.3 to show that there exists a triple that in addition satisfies
$b - a < n$ and $c - b < n$.  (Recall that $a, b, c$ will 
be pairwise noncongruent modulo $n$ by Remark 2.2.)
\endremark

The next result will be useful in the proof of Theorem 2.7.

\proclaim{Lemma 2.5}
\item{\rm (a)}{Let $w \in W$.  We have $$
w(i) < w(i+1) \Leftrightarrow \ell(w s_i) > \ell(w)
$$ and $$
w^{-1}(i) < w^{-1}(i+1) \Leftrightarrow \ell(s_i w) > \ell(w)
.$$}
\item{\rm (b)}{If $c, d \in \zed$ with $c < d$, and $s \in S$ is a generator
of $W$, then $s(c) \geq s(d)$ implies that $c = d-1$, $s(c) = d$ and 
$s(d) = c$.}
\endproclaim

\demo{Proof}
Part (a) is a restatement of \cite{{\bf 7}, Corollary 1.3.3}.  

It is clear from Proposition 1.1 that $s(c) \geq s(d)$ can only occur if
$d - c \leq 2$, because $|s(z) - z| \leq 1$ for all $z \in \zed$.
If $d - c = 1$ and $s(c) \geq s(d)$, then the fact that $s$ induces a 
permutation of $\zed$ implies that $s(c) = d$ and $s(d) = c$ as claimed.
If $d - c = 2$ and $s(c) \geq s(d)$, we must have $s(c) = c + 1$ and $s(d)
= d - 1$, which also cannot happen since $s$ induces a permutation and 
$c + 1 = d - 1$.  Part (b) follows.
\qed\enddemo

The following result is reminiscent of \cite{{\bf 2}, Lemma 1}.

\proclaim{Lemma 2.6}
Let $w \in W$ be fully commutative, and let $s_{i_1} s_{i_2} \cdots
s_{i_r}$ be a reduced expression for $w$.  Suppose $s_{i_j}$ and
$s_{i_k}$ are consecutive occurrences of the generator $s$ (meaning
that $s_{i_j} = s = s_{i_k}$ for some $j < k$ and $s_{i_l} \ne s$ for
all $j < l < k$).  Let $t$ and $t'$ be the two distinct Coxeter
generators that do not commute with $s$.  Then there is an occurrence
$s_{i_p} = t$ and an occurrence $s_{i_q} = t'$ with $j < p, q < k$.
\endproclaim

\demo{Note}
This result is a consequence of \cite{{\bf 5}, Lemma 4.3.5 (ii)} and its
proof, but here we present a more direct argument based on property
R3 of \cite{{\bf 4}, \S2}.
\enddemo

\demo{Proof}
Since $w \in W_c$, we can see that between the two named occurrences of $s$
there must be at least two occurrences of generators that do not
commute with $s$.  If we choose $s$ such that $k - j$ is minimal,
there must be at least one occurrence of each of the generators $t$
and $t'$ in the given interval, as required.
\qed\enddemo

\proclaim{Theorem 2.7}
Let $w \in W = W(\widetilde A_{n-1})$.  The following are equivalent:

\item{\rm (i)}{$w$ is fully commutative (considering $W$ as a Coxeter group);}
\item{\rm (ii)}{if $a, b \in \zed$ with $a < b$ and we have $w(a) > w(b)$,
then $w(a) > a$ and $w(b) < b$;}
\item{\rm (iii)}{$w$ is $321$-avoiding (considering $W$ as a George group);}
\item{\rm (iv)}{there are no positive roots $\a, \be, \a + \be$ 
in the root system of type $\widetilde A_{n-1}$ such that $w(\a) < 0$
and $w(\be) < 0$.}
\endproclaim

\demo{Note}
Recall that the action of a generator $s_i \in W$ on a simple root $\a_j$ is
given by $$
s_i . \a_j = \cases
-\a_j & \text{ if } i = j,\cr
\a_j + \a_i & \text{ if } \bar{i} = \overline{j-1} \text{ or } 
\bar{i} = \overline{j+1},\cr
\a_j & \text{ otherwise}.\cr
\endcases
$$  The reader is referred to \cite{{\bf 6}} or \cite{{\bf 15}} for further
details of the root system.
\enddemo

\demo{Proof}
The equivalence of (i) and (iv) follows from \cite{{\bf 6}, Theorem 2.4}.

\noindent (i) $\Rightarrow$ (ii).

Suppose that (ii) fails to hold.
Suppose $a < b$ and $w(a) > w(b)$ and $w(a) \leq a$.  (The other
possibility, that $w(b) \geq b$, is proved using an argument similar
to the following.)

Let $s_{i_r} s_{i_{r-1}} \cdots s_{i_1}$ be a reduced expression for
$w$.  For $1 \leq j \leq r$, define $$
w[j] := s_{i_j} s_{i_{j-1}} \cdots s_{i_1}
,$$ with $w[0] := 1$ for notational convenience.
Since $a < b$ and $w(a) > w(b)$, there must exist $1 \leq j \leq r$ 
such that $w[j](a) \geq w[j](b)$ and $w[j-1](a) < w[j-1](b)$.  Applying
Lemma 2.5 (b) with $c = w[j-1](a)$, $d = w[j-1](b)$ and $s = s_{i_j}$, we find
that $w[j](b) = w[j-1](a) = w[j-1](b) - 1 = w[j](a) - 1$.  There are
two cases to check; they cover the possibilities because $w(a) \leq
a$.

Case 1: $w(a) < w[j](a)$.  Now $w(a) = w[r](a) < w[j](a) >
w[j-1](a)$.  It follows that there exist $k$ and $l$ with 
$j \leq k < l \leq r$ such that $$
w[k-1](a) + 1 = w[k](a) = w[k+1](a) = \cdots = w[l-2](a) = w[l-1](a) =
w[l](a) + 1
.$$  It follows that $s_{i_k} = s_{i_l}$ and that no generator
$s_{i_p}$ for $k < p < l$ is equal to $s_{i_k}$ or to $s_{1 + i_k}$
(where addition is taken modulo $n$).  Applying Lemma 2.6 to $s_{i_k}$
and $s_{i_l}$ now shows that $w$ is not fully commutative, because
$s_{1 + i_k}$ does not commute with $s_{i_k}$.

Case 2: $a > w[j-1](a)$.  Here $w[0](a) = a > w[j-1](a) < w[j](a)$.
It follows that there exist $k$ and $l$ with $1 \leq k < l \leq j$ and $$
w[k-1](a) - 1 = w[k](a) = w[k+1](a) = \cdots = w[l-2](a) = w[l-1](a) =
w[l](a) - 1
.$$  The consideration of this case is now similar to case 1, mutatis mutandis.

\noindent (ii) $\Rightarrow$ (iii).

Suppose (iii) fails, so that there exist $a < b < c$ with $w(a) > w(b)
> w(c)$.  If $w(b) \leq b$ then (ii) fails when applied to $b < c$,
and if $w(b) \geq b$ then (ii) fails when applied to $a < b$.

\noindent (iii) $\Rightarrow$ (i).

Suppose (i) fails.  Let $s_{i_r} s_{i_{r-1}} \cdots s_{i_1}$ be 
a reduced expression for $w$ for which $s_{i_{t+1}} = s_{i_{t-1}}$ for
some $1 < t < r$.  Define $w[j]$ as in the proof of (i) $\Rightarrow$
(ii) above.  Let us assume that $i_t = i_{t+1} + 1$; the other
case ($i_t = i_{t+1} - 1$) follows by a symmetric argument.

Let $y = w[t-2]$ and $m = i_{t-1}$ (so that $m + 1 = i_t$ by the
previous paragraph).
We first claim that $y^{-1}(m) < y^{-1}(m+1) < y^{-1}(m+2)$.  Since
$s_{i_{t-1}} y > y$, Lemma 2.5 (a) shows that $y^{-1}(m) < y^{-1}(m+1)$.
Since $\ell(s_{i_{t+1}} s_{i_t} s_{i_{t-1}} y) = \ell(y) + 3$ and
$s_{i_{t+1}} s_{i_t} s_{i_{t-1}} = s_{i_t} s_{i_{t-1}} s_{i_t}$, it
follows that $s_{i_t} y > y$ and Lemma 2.5 (a) shows that $y^{-1}(m+1) <
y^{-1}(m+2)$.  (This uses the fact that $i_t = m + 1$.)

Define $a = y^{-1}(m)$, $b = y^{-1}(m+1)$ and $c = y^{-1}(m+2)$, and
let $z = w[t+1]$.  It follows from the definition of the realization of
the generators of $W$ as permutations (given in Proposition 1.1)
and the above paragraph that $z(a) = m+2$, $z(b) = m+1$ and
$z(c) = m$.  We are done once we show that $w(a) > w(b) > w(c)$.
Suppose $w(a) \leq w(b)$.  Then there exists $k$ such that $t + 1 <
k \leq r$ where we have $w[k-1](a) > w[k-1](b)$ and $w[k](a) < w[k](b)$.  This
implies that $w[k-1](a) = w[k-1](b) + 1$ by Lemma 2.5 (b).
Since $\ell(s_{i_k} w[k-1]) > \ell(w[k-1])$, applying Lemma 2.5 (a) 
with $i = w[k-1](b)$ yields a contradiction.  We deduce that $w(a) >
w(b)$ and the proof that $w(b) > w(c)$ is similar.
\qed\enddemo

It is possible to prove the equivalence of conditions (iii) and (iv)
of Theorem 2.7 by using the formalism of George groups from
\cite{{\bf 3}}.  If $w \in W$ then the condition that $w\a < 0$ for a
positive root $\a$ is equivalent to the condition that $\ell(wt) <
\ell(w)$, where $t$ is the (not necessarily simple) reflection
corresponding to the root $\a$.  
The simple reflections are defined in Proposition 1.1 and reflections
in general will be conjugates of these, so they are easily identified.
A more direct approach using roots as opposed to reflections is
implicit in \cite{{\bf 15}, Theorem 1}.

It would be interesting to have an interpretation of condition (ii) of
Theorem 2.7 in terms of roots, since this would enable the condition to be
generalized to other types.  Condition (ii) is also interesting in
that it captures the notion of $321$-avoiding by considering the
images of only two elements.

\head 3. Kazhdan--Lusztig cells \endhead

We wish to apply Theorem 2.7 in order to understand better the
Kazhdan--Lusztig cells in type $\widetilde A_{n-1}$.  For any Coxeter group $W$,
Kazhdan and Lusztig \cite{{\bf 13}} defined partitions of $W$ into left
cells, right cells and two-sided cells, where the two-sided cells are
unions of left (or right) cells.  These cells are naturally
partially ordered by an explicitly defined order, $\leq_{LR}$.  We do
not present the full definitions here; an elementary introduction to
the theory may be found in \cite{{\bf 12}, \S7}.

The properties of these cells are subtle and it seems impossible to
understand them fully without using heavy machinery such as
intersection cohomology.  However, the problem of classifying the
two-sided cells is combinatorially tractable in special cases, and Shi
\cite{{\bf 16}} succeeded in describing combinatorially 
the two-sided cells for the Coxeter group of type $\widetilde A_{n-1}$.
One of the main results of \cite{{\bf 16}} describes a natural bijection between
two-sided cells and partitions of $n$.  This bijection, which is a
generalization of the Robinson--Schensted correspondence, depends on
concepts of which $321$-avoidance is a special case.  

Our aim in \S3 is to explain why the fully commutative elements in
type $\widetilde A$ form a union of Kazhdan--Lusztig cells closed under
the order $\geq_{LR}$.  Fan and Stembridge \cite{{\bf 6}, Theorem 3.1} have
shown that for Coxeter groups of types $A$, $D$, $E$ and $\widetilde A$,
$W_c$ is a union of Spaltenstein--Springer--Steinberg cells.  It is
not always true (even in the aforementioned cases) that $W_c$ is a 
union of Kazhdan--Lusztig cells, but
the situation for finite Coxeter groups is well understood: the
author and J. Losonczy proved in \cite{{\bf 10}, Corollary 3.1.3} that an
irreducible finite Coxeter group has this property if and only if it
does not contain a parabolic subgroup of type $D_4$.

We start by recalling Shi's map $\s$ from $W$ to the partitions of $n$
as described in \cite{{\bf 16}}.  (We have changed the action of $W$ on
$\zed$ from being a right action to being a left action, but this
makes no difference and $\s(w) = \s(w^{-1})$.)

\definition{Definition 3.1}
Fix $w \in W = W(\widetilde A_{n-1})$.  Define, for $k \geq 1$, $$
d_k  = d_k(w) := \max \left\{ |X| : X = \bigcup_{i = 1}^k X_i \subset
\zed\right\}
$$ where (a) the union is disjoint,
(b) any two integers $u, v \in X$ are required not to be
congruent modulo $n$, and (c) whenever the subset $X_i$ of $\zed$ contains 
integers $u$ and
$v$ with $u < v$ then we require $w(u) > w(v)$.  (Here,
$|X|$ denotes the cardinality of $X$.)

The partition $\s(w)$ of $n$ is defined by $$
(d_1, d_2 - d_1, d_3 - d_2, \ldots, d_t - d_{t-1})
,$$ where $t$ is maximal subject to the condition that $d_{t-1} < d_t$
(meaning that $d_t = n$).
\enddefinition

It may not be obvious from the definition that $\s(w)$ is a partition
(\idest $d_1 \geq d_2 - d_1 \geq \cdots \geq d_t - d_{t-1}$), but this
follows by a result of Greene and Kleitman \cite{{\bf 11}}.

Recall that if $\l$ and $\mu$ are two partitions of $n$, we say that
$\l \domeq \mu$ ($\l$ dominates $\mu$) if for all $k$ we have $$
\sum_{i = 1}^k \l_i \geq \sum_{i = 1}^k \mu_i
.$$  This defines a partial order on the set of partitions of $n$.
As mentioned before, there is also a natural partial order $\leq_{LR}$
on the two-sided Kazhdan--Lusztig cells of a Coxeter group.  The
following theorem explains how these two orders are related.

\proclaim{Theorem 3.2 (Shi)}
Let $y, w \in W = W(\widetilde A_{n-1})$.  Then $y \leq_{LR} w$ in the
sense of Kazhdan--Lusztig if and only if $\s(y) \domeq \s(w)$.  In
particular, $y$ and $w$ are in the same two-sided cell (\idest $y
\leq_{LR} w \leq_{LR} y$) if and only if $\s(y) = \s(w)$.
\endproclaim

\demo{Proof}
See \cite{{\bf 17}, \S2.9}.
\qed\enddemo

As we now explain, it is not hard to see from these definitions that
the $321$-avoiding elements of $W$ are a union of Kazhdan--Lusztig
cells.

\proclaim{Lemma 3.3}
Fix $w \in W$.
\item{\rm (i)}{If $w$ is $321$-avoiding then all the sets $X_i$ in
Definition 3.1 satisfy $|X_i| \leq 2$.}
\item{\rm (ii)}{If $w$ is not $321$-avoiding then the number $d_1$ in
Definition 3.1 satisfies $d_1 \geq 3$.}
\endproclaim

\demo{Proof}
If there is a set $X_i$ with $|X_i| \geq 3$ then by definition of
$X_i$, there exist $u_1, u_2, u_3 \in X_i$ with $u_1 < u_2 < u_3$ 
but $w(u_1) > w(u_2) > w(u_3)$.  This
proves (i).  The proof of (ii) is similar: if $a < b < c$ with $w(a) >
w(b) > w(c)$ then we may choose $X_1$ to contain these integers $a, b,
c$ (and possibly others), thus showing that $d_1 \geq 3$.
\qed\enddemo

\proclaim{Theorem 3.4}
Let $W_c$ denote the set of fully commutative elements of $W(\widetilde
A_{n-1})$.  Then $W_c$ is a union of two-sided Kazhdan--Lusztig cells
closed under $\geq_{LR}$.
\endproclaim

\demo{Note}
By the statement $W_c$ is closed under $\geq_{LR}$ we mean that if $w
\in W_c$ and $y \geq_{LR} w$ then $y \in W_c$.
\enddemo

\demo{Proof}
It is clearly enough to show that if $w \in W\backslash W_c$ and 
$w \geq_{LR} y$ then $y \in W \backslash W_c$.  By theorems 2.7 and
3.2, it is enough to show that if $w$ is not $321$-avoiding
and $\s(y) \domeq \s(w)$ then $y$ is not $321$-avoiding.

Suppose $w$ and $y$ are as above. 
Let $\l = \s(w)$ and $\mu = \s(y)$.  By Lemma 3.3 (ii), $\l_1 \geq 3$.
Since $\mu \domeq \l$, it follows from the definition of the dominance
order that $\mu_1 \geq 3$.  By Lemma 3.3 (i), this means that $y$ is
not $321$-avoiding, as required.
\qed\enddemo

We will discuss some of the applications of Theorem 3.4 in \S5.

\head 4. The extended affine Weyl group associated to $GL_n(\complex)$ \endhead

We now return to the group $\W$ of Proposition
1.3.  The generalization of Theorem 2.7 to this group is
straightforward as we now show.  Recall that every element $w' \in \W$
may be written uniquely as $\rho^z w$ for $z \in \zed$ and $w \in W$,
where $\rho$ is as in Definition 1.2.

\definition{Definition 4.1}
Let $w' = \rho^z w$ be a typical element of $\W$.  We say $w'$ is {\it
fully commutative} if $w$ is fully commutative, and we say $w'$ is
{\it $321$-avoiding} if $w$ is $321$-avoiding.
\enddefinition

The definition of $321$-avoiding above is a natural one as if $a < b$
then $w(a) < w(b)$ if and only if $w'(a) < w'(b)$.

The group $\W$ may be made to act on the root system of $W$ in a
natural way which is clear from the following fact.

\proclaim{Lemma 4.2}
Let $s_i$ be a Coxeter generator of $W$.  Then $\rho s_i \rho^{-1} =
s_{i+1}$, where addition is taken modulo $n$.
\endproclaim

\demo{Proof}
This follows by considering the induced permutations on
$\zed$.
\qed\enddemo

We define the action of $\W$ on the root system by stipulating that 
$\rho(\a_i) = \a_{i+1}$ (with addition taken modulo $n$); this is 
reasonable by Lemma 4.2.  Notice that
$\rho^n$ acts as the identity; the extended affine Weyl group
associated to $SL_n(\complex)$ is the
quotient of $\W$ by the relation $\rho^n = 1$.

\proclaim{Proposition 4.3}
Let $w \in \W$.  The following are equivalent:

\item{\rm (i)}{$w$ is fully commutative;}
\item{\rm (ii)}{$w$ is $321$-avoiding;}
\item{\rm (iii)}{there are no positive roots $\a, \be, \a + \be$ 
in the root system of type $\widetilde A_{n-1}$ such that $w(\a) < 0$
and $w(\be) < 0$.}
\endproclaim

\demo{Proof}
The equivalence of (i) and (ii) is immediate from the definitions and
Theorem 2.7.  It follows from the definition of the action of $\rho$
on simple roots that powers of $\rho$ permute the positive roots linearly.
This means that if $w = \rho^z w'$ then $w'(\a) < 0$ and $w'(\be) < 0$ if
and only if $w(\a) < 0$ and $w(\be) < 0$.  The equivalence of (iii) with
the other conditions now follows by Theorem 2.7.
\qed\enddemo

There seems to be no direct analogue of condition (ii) of Theorem 2.7
for the group $\W$.

Our results here suggest the following question about George groups.

\proclaim{Question 4.4}
Is it possible to define in a uniform way an extended George group for
affine types $A$, $B$, $C$ and $D$ that specializes to the group $\W$ in 
type $A$?
\endproclaim

We do not necessarily expect that the extended George group for types $B$, 
$C$ and $D$ should be closely related to the group $\W$.

It is also natural to wonder whether Theorem 3.4 may be extended to
other affine types.  The theorem will fail for Coxeter systems
containing a parabolic subgroup of type $D_4$ by the results of
\cite{{\bf 10}}, which rules out type $\widetilde D$ and type $\widetilde B_l$
(at least for $l \geq 4$).  The most interesting open case is
therefore that of type $\widetilde C$.

\head 5. Applications to Kazhdan--Lusztig theory \endhead

Theorem 3.4 shows that for the group $W = W(\widetilde A_{n-1})$, the set
$W_c$ is a union of two-sided Kazhdan--Lusztig cells.  This result may be
refined as follows.

\proclaim{Proposition 5.1}
Let $W_c$ denote the set of fully commutative elements of \newline 
$W(\widetilde
A_{n-1})$.  Then the number of distinct two-sided cells contained in $W_c$ is
$(n+1)/2$ if $n$ is odd, and $(n+2)/2$ if $n$ is even.  A set of 
representatives in $W$ for the two-sided cells is given by the set $$
\left\{ s_2 s_4 \cdots s_{2k} : 0 \leq k \leq {n \over 2} \right\}
,$$ using the usual numbering of the generators shown in Figure 1.
\endproclaim

\demo{Proof}
By theorems 2.7, 3.2, 3.4 and Lemma 3.3, the number of fully commutative 
two-sided cells in $W_c$ is equal to the number of partitions of $n$ with
all parts less than or equal to $2$, and the first claim follows.  It follows
from Definition 3.1 that if $w = s_2 s_4 \cdots s_{2k}$, the partition $\s(w)$
has $k$ parts equal to $2$ and the other parts equal to $1$.  The result
follows.
\qed\enddemo

We end by discussing briefly the application of Theorem 3.4 to the computation
of certain structure constants for the Kazhdan--Lusztig basis.

The Kazhdan--Lusztig
basis, which first appeared in \cite{{\bf 13}}, is a free 
$\zed[v, v^{-1}]$-basis for the Hecke algebra $\H(W)$ associated to the 
group $W$.  The basis, $\{C'_w : w \in W\}$, is naturally indexed by $W$.
The structure constants $g_{x, y, z}$, namely the Laurent polynomials
occurring in the expression $$
C'_x C'_y = \sum_{z \in W} g_{x, y, z} C'_z
,$$ have many subtle properties, such as the fact that $g_{x, y, z} \in
\enn[v, v^{-1}]$.  (It seems impossible to establish this fact 
combinatorially.)

Theorem 3.4 implies that the $\zed[v, v^{-1}]$-span of the set $$
J(W) = \{ C'_w : w \not\in W_c \}
$$ is an ideal of $\H(W)$.  The quotient algebra $\H(W)/J(W)$
is equipped with a natural basis $\{c_w : w \in W_c\}$, where
$c_w := C'_w + J(W)$.  The basis $\{c_w : w \in W_c\}$ is the canonical basis 
(in the sense of \cite{{\bf 9}}) for the quotient of the Hecke algebra $\H(W)/J(W)$
by \cite{{\bf 10}, Theorem 2.2.3}.  Although we do not present full details,
it is possible to realize this canonical basis combinatorially; the
description is similar to that given in \cite{{\bf 8}, Definition 6.4.3}.  This
means that, if $z \in W_c$, we can compute 
$g_{x, y, z}$ by simple
combinatorial means.  (Note that if $z \in W_c$ and $g_{x, y, z} \ne 0$ then
we must have $x, y \in W_c$.)  In particular, it is easily checked that
$g_{x, y, z} \in \enn[v, v^{-1}]$ in this case.

\head Acknowledgements \endhead

The author thanks J.R. Stembridge for some helpful discussions and
J. Losonczy for comments on the manuscript.

\leftheadtext{}
\rightheadtext{}

\end